\documentclass[a4paper,11pt]{article}
\usepackage{enumerate}
\usepackage{latexsym}
\usepackage{amsfonts}
\usepackage[only,ninrm,elvrm,twlrm,sixrm,egtrm,tenrm]{rawfonts}
\usepackage{indentfirst}
\usepackage{amsmath}
\usepackage[noend]{algorithmic}
\usepackage{algorithm}
\usepackage{graphicx,psfrag}
\usepackage{graphics}
\usepackage{makeidx}
\newtheorem{thm}{Theorem}[section]

\newtheorem{lem}[thm]{Lemma}

\newtheorem{prop}[thm]{Proposition}

\def\qed{\hfill \nopagebreak\rule{5pt}{8pt}}
\title{\bf Partitioning complete graphs by heterochromatic trees
 \footnote{Supported by NSFC, PCSIRT and the ``973" program. }}

\author{
\small  Zemin Jin$^1$ and Xueliang Li$^2$ \\
[3mm] \small    $^{1}$Department of Mathematics, Zhejiang Normal University\\
\small Jinhua 321004, P.R. China\\
\small $^{2}$Center for Combinatorics and LPMC,  Nankai University\\
\small Tianjin 300071, P.R. China\\}

\begin{document}

\maketitle
\begin{abstract}

A {\it heterochromatic tree} is an edge-colored tree in which any
two edges have different colors. The {\it heterochromatic tree
partition number} of an $r$-edge-colored graph $G$, denoted by
$t_r(G)$, is the minimum positive integer $p$ such that whenever
the edges of the graph $G$ are colored with $r$ colors, the
vertices of $G$ can be covered by at most $p$ vertex-disjoint
heterochromatic trees. In this paper we determine the
heterochromatic tree partition number of an $r$-edge-colored
complete graph.\\
[0.1in] {\bf  Keywords:}  edge-colored graph, heterochromatic
tree, partition.\\
[2mm] {\bf AMS subject classification (2000)}: 05C05, 05C15, 05C70,
05C75.
\end{abstract}

\section{Introduction}
A {\it monochromatic (heterochromatic) tree} is an edge-colored
tree in which any two edges have the same (different) color(s).
The {\it (monochromatic) tree partition number} of an
$r$-edge-colored graph $G$ is defined to be the minimum positive
integer $p$ such that whenever the edges of $G$ are colored with
$r$ colors, the vertices of $G$ can be covered by at most $p$
vertex-disjoint monochromatic trees. The {\it (monochromatic)
cycle partition number} and the {\it (monochromatic)
path partition number} are defined similarly.\\

Erd\H{o}s, Gy\'{a}rf\'{a}s and Pyber \cite{egp} proved that the
(monochromatic) cycle partition number of an $r$-edge-colored
complete graph $K_n$ is at most $cr^2\ln r $ for some constant
$c$. This implies a conjecture from \cite{gys} in a stronger form.
Recently, the bound was improved by Gy\'{a}rf\'{a}s et al.
\cite{gyfs}. Almost solving one of the conjectures in \cite{egp},
Haxell and Kohayakawa \cite{hk} proved that the (monochromatic)
tree partition number of an $r$-edge-colored complete graph $K_n$
is at most $r$ provided that $n$ is large enough with respect to
$r$. Haxell \cite{ha} proved that the (monochromatic) cycle
partition number of an $r$-edge-colored complete bipartite graph
$K_{n,n}$ is also independent of $n$, which answered a question in \cite{egp}.\\

From above, one can see that the (monochromatic) tree, path, and
cycle partition number of $r$-edge-colored graphs $K_{n}$ and
$K_{n,n}$ are independent of $n$. The same seems to be not true
for other graphs. Also, no (monochromatic) partition number of an
$r$-edge-colored graph $K_{n}$ or $K_{n,n}$ is determined exactly.
The only exception is due to Kaneko, Kano and Suzuki \cite{kks},
who gave an explicit expression for the (monochromatic) tree
partition number of a 2-edge-colored complete multipartite graph.
In particular, let $n_1, n_2, \cdots, n_k$ ($k\geq 2$) be integers
such that $1\leq n_1\leq n_2\leq \cdots \leq n_k$ and let $n=n_1+
n_2+\cdots +n_{k-1}, m=n_k$. The authors \cite{kks} proved that
$$
t_2^{'}(K_{n_1, n_2, \cdots, n_k})=\lfloor \frac{m-2}{2^n}
\rfloor+2,
$$
where $t_r^{'}(K_{n_1, n_2, \cdots, n_k})$ denotes the
(monochromatic) tree partition number of the $r$-edge-colored graph
$K_{n_1, n_2, \cdots, n_k}$. Other related partition
problems can be found in \cite{eg,Luczak,rado}.\\

Analogous to the monochromatic tree partition case, the authors
\cite{chen} introduced the definition of {\it heterochromatic tree
partition number} of an $r$-edge-colored graph $G$. The {\it
heterochromatic tree partition number} of an $r$-edge-colored
graph $G$, denoted by $t_r(G)$, is defined to be the minimum
positive integer $p$ such that whenever the edges of the graph $G$
are colored with $r$ colors, the vertices of $G$ can be covered by
at most $p$ vertex-disjoint heterochromatic trees. In \cite{chen},
the authors determined the heterochromatic tree partition number
of an $r$-edge-colored complete bipartite graph $K_{m,n}$. In this
paper we consider an $r$-edge-colored complete graph $K_{n}$ and give
the exact expression for its heterochromatic tree partition number. \\

Before proceeding, we introduce some definitions and notations.
Throughout this paper, we use $r$ to denote the number of the
colors, and an $r$-edge-coloring of a graph $G$ means that each
color appears at least once in $G$. Let $\phi$ be an
$r$-edge-coloring of a graph $G$. For an edge $e\in E(G)$, denote
by $\phi(e)$ the color of $e$. Denote by $t_r(G, \phi)$ the
minimum positive integer $p$ such that under the $r$-edge-coloring
$\phi$, the vertices of $G$ can be covered by at most $p$
vertex-disjoint heterochromatic trees. Clearly, $t_r(G)=\max
_{\phi} t_r(G, \phi)$, where $\phi$ runs over all
$r$-edge-colorings of the graph $G$. Let $\phi$ be an
$r$-edge-coloring of the graph $G$ and $F$ be a spanning forest of
$G$, each component of which is a heterochromatic tree. If $F$
contains exactly $t_r(G, \phi)$ components, then $F$ is called an
{\it optimal heterochromatic tree partition} of the graph $G$ with
edge-coloring $\phi$. Note that a tree consisting of a single
vertex is also regarded as a heterochromatic tree. \\

For any integer $r\geq 2$, there is a unique positive integer $t$,
such that ${t\choose 2} +2\leq r\leq {{t+1}\choose 2} +1$.
Clearly, the integer $t$ is determined completely by $r$, and here
we denote it by $f(r)=t$. This integer $f(r)=t$ will play an
important role in expressing the number $t_r(K_n)$. If the color
number $r=1$, clearly a maximum matching (plus a single vertex
when $n$ is odd) in $K_n$ is an optimal heterochromatic tree
partition, and then $t_r(K_n)=\lceil \frac{n}{2} \rceil$. So, in
the rest of this paper we only consider the case $2\leq r\leq
{n\choose 2}$. The following is the main result of this paper.

\begin{thm}\label{main}
Let $n\geq 3$, $2\leq r\leq {n\choose 2}$ and $f(r)=t$. Then
$t_r(K_n)=\lceil \frac{n-t}{2} \rceil$.
\end{thm}

As we know, the monochromatic tree partition number of an
edge-colored complete graph $K_n$ is bounded by a function
independent of $n$, and from the result mentioned above, the
heterochromatic tree partition number does not have this property
any more. The rest of the paper is organized as follows. In
Section 2, we present a canonical $r$-edge-coloring $\phi^{*}$ of
the complete graph $K_{n}$, and show that, under the canonical
$r$-edge-coloring $\phi^{*}$, the optimal heterochromatic tree
partition in the graph $K_{n}$ contains exactly $\lceil
\frac{n-t}{2} \rceil$ components. The proof of our main result is
complete in the last section.

\section{A canonical $r$-edge-coloring $\phi^{*}_r$}

In this section we present a canonical $r$-edge-coloring of the
graph $K_n$. Let $f(r)=t$, i.e., ${t\choose 2} +2\leq r\leq
{{t+1}\choose 2} +1$. Let $S\subseteq V(K_n)$ and $|S|=t$. Take a
vertex $u\in V(K_n)-S$. We define the canonical $r$-edge-coloring
$\phi^{*}_r$ by
\begin{enumerate}
\item  giving distinct colors to the edges of $K_n[S]$; \item for
each color not used, assign it to an edge $uv$ (if it is not
colored) between $u$ and $S$; \item finally, color all the
remaining edges by the color not used if it exists, or else by the
same color.
\end{enumerate}

We have the following proposition.

\begin{prop}\label{prop}
$t_r(K_n, \phi^{*}_r)=\lceil \frac{n-t}{2}\rceil.$
\end{prop}
{\it Proof:} First, we present an heterochromatic tree partition
with exact $\lceil \frac{n-t}{2}\rceil$ components, which implies
that $t_r(K_n, \phi^{*}_r)\leq \lceil \frac{n-t}{2} \rceil$. Let
$X=S\cup \{u\}\cup \{v\}$, where $v\in V(K_n)-S-u$. It is easy to
see that $K_n[X]$ contains a heterochromatic spanning tree $T$,
and the vertices not in $X$ induce a monochromatic complete
subgraph which can be covered by $\lceil \frac{n-t-2}{2}\rceil$
disjoint heterochromatic trees. So, the union of $T$ and those
$\lceil \frac{n-t-2}{2}\rceil$ disjoint heterochromatic trees
consist of a heterochromatic tree partition of $K_n$. This implies
that $t_r(K_n, \phi^{*}_r)\leq \lceil \frac{n-t}{2} \rceil$.\\

Next, we prove that $t_r(K_n, \phi^{*}_r)\geq \lceil \frac{n-t}{2}
\rceil$. Suppose on the contrary that $t_r(K_n, \phi^{*}_r)<
\lceil \frac{n-t}{2}\rceil$ for some $n$ and $r$. \\

Let $F$ be an optimal heterochromatic tree partition of $K_n$ with
$r$-edge-coloring $\phi^{*}_r$. Denote by $T_1, T_2, \cdots, T_k$
the components of $F$ which contains vertices of $S$. We choose
$F$ such that the number of trees covering $S$ is as small as
possible. Note that each component of $F$ not containing any
vertex of $S$ is an edge or a single vertex, and at most one of
the components of $F$ is a single vertex. Since $F$ is an optimal
heterochromatic tree partition, from the definition of
$\phi^{*}_r$, we have the following facts.\\
{\bf Fact 1.} $u\in T_i$ for some $1\leq i\leq k$.\\
{\bf Fact 2.} $|T_j\cap (V(K_n)-S-u)|=1$ for each $T_j$.\\

If $k=1$, it is easy to see that $F$ just contains $\lceil
\frac{n-t}{2} \rceil $ trees.  So, assume that $k\geq 2$. Let
$S\cap T_i=S_i$ and $v_i=T_i\cap (V(K_n)-S-u)$. From the
definition of $\phi^{*}_r$, we have that there exists a
heterochromatic tree, denoted by $T$, covering all the vertices
$T_1\cup (T_2-v_2)$. So $F^{'}=(F-T_1-T_2)\cup \{T\} \cup \{v_2\}$
is an optimal heterochromatic tree partition such that the number
of trees covering $S$ is $k-1$, a contradiction, which completes the
proof. \qed \\

\section{Proof of Theorem  \ref{main}}

Given a complete graph $K_n$, the heterochromatic tree partition
number is closely related to the color number. Before proving our
main result, we have the following lemma which presents the
relationship between $t_{r+1}(K_n)$ and $ t_{r}(K_n)$.

\begin{lem} \label{inequ}
$t_{r+1}(K_n)\leq t_{r}(K_n)$.
\end{lem}
{\it Proof:} Given any $(r+1)$-edge-coloring $\varphi$ of $K_n$.
Denote by $E_{i}$ the set of edges colored by the color $i$.
Recoloring the edges of $E_{r+1}$ by the color $r$, we obtain a
$r$-edge-coloring $\psi$ of $K_n $. Clearly, $t_{r+1}(K_n,
\varphi)\leq t_{r}(K_n, \psi)$. So, $t_{r+1}(K_n)\leq t_{r}(K_n)$. \qed \\

The following lemma presents the relationship between the
edge-connectivity and size of a graph. The proof is omitted here.
\begin{lem}\label{cutedge}
Let $G$ be a simple graph of order $n$. If $G$ contains a
cut-edge, then $|E(G)|\leq {{n-1}\choose 2}+1$.\qed \\
\end{lem}

\noindent {\it  Proof of Theorem  \ref{main}:} We prove the
theorem by induction on $r$ and $n$. First, we consider the case
$r=2$. Let $\phi$ be a $2$-edge-coloring of $K_n$. Note that for
any $2$-edge-coloring of $K_n$, $n\geq 3$, there is always a
heterochromatic tree of order three. Then, we can easily find
$1+\lceil \frac{n-3}{2}\rceil=\lceil \frac{n-1}{2}\rceil $
vertex-disjoint heterochromatic trees which cover all the
vertices. So we have $t_r(K_n, \phi)\leq \lceil
\frac{n-1}{2}\rceil $. Then, from Proposition \ref{prop} the
result holds for $r=2$. Obviously, the result holds for $n=3,4$. \\

Assume that the result holds for the color number less than $r$ or
the order of a complete graph less than $n$. Now we consider the
$r$-edge-colored complete graph $K_n$, $r\geq 3$. Let $f(r)=t$. If
${t\choose 2} +3\leq r\leq {t+1\choose 2} +1$, then $f(r-1)=t$. By
the induction hypothesis, $t_{r-1}(K_n)=\lceil
\frac{n-t}{2}\rceil$. From Lemma \ref{inequ}, $t_r(K_n)\leq
t_{r-1}(K_n)=\lceil \frac{n-t}{2}\rceil$. And, from Proposition
\ref{prop}, $t_r(K_n)\geq t_r(K_n, \phi^{*}_r)=\lceil
\frac{n-t}{2}\rceil$. Then, we have
$t_r(K_n)=\lceil \frac{n-t}{2}\rceil$, as desired.\\

So, we only need to consider the case $r={t\choose 2} +2$. Let
$\phi$ be an $r$-edge-coloring of $K_n$. Let $G$ be a
heterochromatic subgraph of $K_n$, such that $\delta(G)\geq 1$
and, for each color $i$, there is a unique edge colored by the
color $i$ in $G$. Denote by $G_1, G_2, \cdots, G_k$ the components
of $G$, where the order of $G_i$ is $n_i$, $1\leq i\leq k$, and
$n_1\geq n_2\geq \cdots \geq n_k\geq 2$. Choose $G$ such that
$n_1$ is as large as possible. Since the color number
$r\geq 3$, we have $n_1\geq 3$.\\

Suppose that $k=1$. By $r={t\choose 2}+2$, we have $n_1\geq t+1$.
If $n_1\geq t+2$, then $t_r(K_n, \phi)\leq 1+\lceil
\frac{n-n_1}{2} \rceil\leq \lceil \frac{n-t}{2}\rceil$. So, assume
$n_1=t+1$. By Lemma \ref{cutedge}, $G$ does not contain any
cut-edge. Let $g\in [V(G_1), \overline{V(G_1)}]$, i.e., one
end-vertex of $g$ belongs to $V(G_1)$ and the other one belongs to
$\overline{V(G_1)}]$.  From the choice of $G$, there is an edge
$h\in E(G_1)$ with $\phi(h)=\phi(g)$. Since $G$ does not contain
any cut-edge, by deleting the edge $h$ and adding the edge $g$, we
can find a heterochromatic graph with $r$ edges, the largest
component of which has an order at least $n_1+1$. A contradiction
to the choice of the graph $G$. \\

So, assume $k\geq 2$. If $n_1\geq t+2$, then $t_r(K_n, \phi)\leq
1+\lceil \frac{n-n_1}{2} \rceil\leq \lceil \frac{n-t}{2}\rceil$,
as desired. Thus, assume $n_1\leq t+1$. We have the following claim.\\

\noindent {\bf Claim: } $G_1$ contains a cut-edge, and then
$|E(G_1)|\leq {n_1-1\choose 2}+1$.

Otherwise, suppose that $G_1$ does not contain any cut-edge. Let
$g\in [V(G_1), \overline{V(G_1)}]$. From the choice of $G$, there
is an edge $h\in E(G_1)$ with $\phi(h)=\phi(g)$. Since $G$ does
not contain any cut-edge, by deleting the edge $h$ and adding the
edge $g$, we can find a heterochromatic graph with $r$ edges, the
largest component of which has an order at least $n_1+1$, a
contradiction to the choice of the graph $G$. From Lemma
\ref{cutedge}, $|E(G_1)|\leq {n_1-1\choose 2}+1$ follows clearly.
This completes the proof  of the claim. \\

Now we consider the graph $K_n-V(G_1)$, a complete graph of order
$n-n_1$. When restricting the $r$-edge-coloring $\phi$ on the
graph $K_n-V(G_1)$, we have that $K_n-V(G_1)$ is edge-colored by
$r_0$ colors, where $r_0\geq r-({n_1-1\choose 2}+1)$. If $r_0\geq
2$, let $f(r_0)=t_0$. It follows that either $r_0=1$, or $t_0\geq
t-n_1+1$. We distinguish the following cases.\\

{\it  Case 1.} $r_0=1$.

Then $K_n-V(G_1)$ is monochromatic, and then it follows that $k=2$
and $n_2=2$. Let $G_2=uv$. From the choice of $G$, we have
$|E(G_1)|=r-1={t\choose 2}+1$. By $n_1\leq t+1$, we have
$n_1=t+1$. From Claim 1, let $e$ be a cut-edge in $G_1$. Since
$|E(G_1)|={t\choose 2}+1$ and $n_1=t+1$, we have $G_1-e\cong
K_{t}\cup K_1$. Let $w\in V(G_1)$. From the choice of $G$, we have
$\phi(uw)\neq \phi(uv)$, and there is a cut-edge in $ G_1$ colored by
the same color $\phi(uw)$.\\

If $n_1\geq 4$, from $G_1-e\cong K_{t}\cup K_1$, we have that $e$
is the unique cut-edge in $G_1$. By $G_1-e\cong K_t\cup K_1$, we
can take a vertex $w$ which is not single in $G_1-e$. Then
$\phi(uw)=\phi(e)$. By deleting the edge $e$ and adding the edge
$uw$, we can find a heterochromatic graph with $r$ edges, the
largest component of which has an order at least $n_1+1$,
a contradiction to the choice of $G$. \\

So, assume $n_1=3$. Then $r=3$ and $G_1\cong P_3$. Let $G_1=xyz$.
Then either $\phi(yu)=\phi(xy)$ or $\phi(yu)=\phi(yz)$. Without
loss of generality, assume $\phi(yu)=\phi(xy)$. Then $\phi(yu)\neq
\phi(yz)$. Again, the graph $zyuv$ is heterochromatic and of size $r$,
a contradiction to the choice of $G$.\\

{\it Case 2.} $t_0\geq t-n_1+1$.

Since $r_0\geq 2$, we have $t_0\geq 1$. If $t_0\geq t-n_1+2$, then
by the induction hypothesis, the graph $K_n-V(G_1)$ can be covered
by at most $\lceil \frac{n-n_1-t_0}{2}\rceil$ vertex-disjoint
heterochromatic trees. Thus, $t_r(K_n, \phi)\leq 1+ \lceil
\frac{n-n_1-t_0}{2}\rceil\leq \lceil \frac{n-t}{2} \rceil $, as desired.\\

Suppose $t_0= t-n_1+1$. Then we have $r={t\choose 2}+2\leq
|E(G_1)|+r_0\leq {n_1-1\choose 2}+1+{t_0+1\choose
2}+1={n_1-1\choose 2}+1+{t-n_1+1+1\choose 2}+1$. This implies that
${t\choose 2}\leq {n_1-1\choose 2}+{t-(n_1-1)+1\choose 2}$, i.e.,
$(n_1-1)(t-(n_1-1))\leq t-(n_1-1)$. By $n_1\geq 3$ and $n_1\leq
t+1$, we have $n_1=t+1$, and then $t_0=0$, a contradiction to the
fact $t_0\geq 1$. The proof is now complete. \qed

\end{document}